\documentclass[a4paper, 12pt]{amsart}
\usepackage{amsmath, amssymb, eucal, amscd, amstext, enumerate}
\usepackage{amsfonts,amsthm}
\usepackage{mathrsfs}

\theoremstyle{plain}
\newtheorem{thm}{Theorem}[section]
\newtheorem{lem}[thm]{Lemma}

\newtheorem{prop}[thm]{Proposition}

\newtheorem{defn}[thm]{Definition}

\newtheorem{rem-ntn}[thm]{Remark and Notation}

\newenvironment{prf}{{\noindent \textbf{Proof:}\ }}{\hfill $\Box$\\ \smallskip}

\numberwithin{equation}{section}

\newcommand{\abs}[1]{\left\vert#1\right\vert}
\newcommand{\norm}[1]{\left\|#1\right\|}

\newcommand{\CL}{\mathcal{L}}

\newcommand{\BR}{\mathbb{R}}

\begin{document}

\title[]{Uniform local Lipschitz continuity of eigenvalues with respect to the potential in $L^1[a,b]$}

\author{Xiao Chen \and Jiangang Qi}
\address{School of Mathematics and Statistics, Shandong University at
         Weihai \\Weihai 264209, P.R. China}
\email{chenxiao@sdu.edu.cn; qjg816@163.com}\




\begin{abstract}
The present paper shows that the eigenvalue sequence $\{\lambda_n(q)\}_{n\geqslant 1}$ of regular Sturm-Liouville
eigenvalue problem with certain monotonic weights  is uniformly Lipschitz continuous with respect to the potential $q$ on any bounded subset of $L^1([a,b],\mathbb{R})$.
\vspace{04pt}

\vspace{04pt}

\end{abstract}

\begin{thanks} {The first named author is supported by the NSF of China (Grant 11701327) and China Postdoctoral Science Foundation (Grant 2017M612252).
The second named author, as the corresponding author, is supported by the NSF of China (Grant 11271229).}
\end{thanks}

\maketitle

\section{Introduction}
\bigskip

Consider the regular Sturm-Liouville eigenvalue problem associated to the second order differential equation
\begin{equation}\label{def:S-L}
-(p(x)y(x)')'+q(x)y(x)=\lambda\omega(x) y(x) \text{ on } [a,b]
\end{equation}
with the self-adjoint separated boundary conditions
\begin{equation}\label{def:S-L-bdry}
y(a)\cos\alpha+(py')(a)\sin\alpha=0,\ y(b)\cos\beta+(py')(b)\sin\beta=0,
\end{equation}
where $\alpha,\ \beta\in[0,\pi),$ $\lambda$  is the spectral parameter,
\begin{equation}\label{condition}
\frac{1}{p},\ q,\ \omega\in L^1([a,b],\mathbb{R}),\quad p,\ \omega>0\ \ a.e.\text{ on }[a,b].
\end{equation}
Here $L^1[a,b]$ denotes the Banach space of all Lesbegue integrable, complex valued functions on the closed interval $[a,b]\subset\BR$
 equipped with the canonical $L^1$-norm $\norm{\cdot}_{L^1}$. The subspace of real valued functions of $L^1[a,b]$ is denoted by $L^1([a,b],\mathbb{R})$.  

Under the natural condition \eqref{condition}, the eigenvalue problem, \eqref{def:S-L}
and \eqref{def:S-L-bdry}, admits only countably infinite  number of real eigenvalues
which are isolated,  bounded below and unbounded above by the spectral theory of differential operators.

Fix $p$ and $\omega$, let $\lambda_n(q)$ be the $n$th eigenvalue with respect to the potential
function $q$. It is well known that
\begin{equation}\label{eigen-seq}
-\infty<\lambda_1(q)<\lambda_2(q)<\cdots<\lambda_n(q)<\cdots,
\end{equation}
and
\begin{equation}\label{eigen-infty}
\lambda_n(q)\to\infty \text{ as } n\to\infty.
\end{equation}
Moreover, $\lambda_n(q)$ can be viewed as a functional on $L^1([a,b],\mathbb{R})$ for every $n\ge 1$.
It is also known that $\lambda_n(q)$ is continuous, and even differentiable,
with respect to $q$ in $L^1[a,b]$ (see e.g. \cite{Zet05} as well as \cite{KZ96}, \cite{MZ96} and \cite{Pos87}).

The continuity and differentiability of eigenvalues provide efficient tools in the
study of properties of eigenvalues and eigenfunctions as well as in other related fields.
In the recent years, Professor Meirong Zhang and his collaborators have obtained fruitful results
on weak and strong continuity of eigenvalues and eigenvalue-pairs of several kinds of
eigenvalue problems
 (see e.g. \cite{ZWMQX18}, \cite{WYZ18}, \cite{CMZ18}, \cite{Zhang14}, \cite{YZ11} and \cite{Zhang08} as well as \cite{MeZ09} and \cite{MeZ10}).

The main topic of this paper is the study of a new continuity, called \emph{uniform local Lipschitz continuity}, 
of the eigenvalue sequence $\{\lambda_n(q)\}_{n\geqslant 1}$ with respect to the potential 
function $q$ in $L^1([a,b],\mathbb{R})$.

\begin{defn}\label{def:subunif-cont-q-SL} The eigenvalue sequence $\{\lambda_n(q)\}_{n\geqslant 1}$ of 
(\ref{def:S-L})-(\ref{def:S-L-bdry}) is said to be \emph{uniformly locally Lipschitz continuous with respect to the potential $q$ in $L^1([a,b],\mathbb{R})$}, if, for any $L^1$-norm bounded subset $\Omega\subset L^1([a,b],\mathbb{R})$,
  there exists a positive number $C(\Omega)$ such that
\begin{equation}\label{def:uni-cont-eigv}
\abs{\lambda_n(q_1)-\lambda_n(q_2)}\leqslant C(\Omega)\norm{q_1-q_2}_{L^1}
\end{equation}
for all $n\geqslant 1$, $q_1,q_2\in\Omega$.
\end{defn}

Note that $C(\Omega)$ is independent of the index $n$ of the eigenvalues $\{\lambda_n(q)\}_{n\geqslant 1}$, and hence this local Lipschitz continuity is uniform for all $n\geqslant 1$. This is exactly the meaning of the word ``uniformly"  in the definition above.

The present paper shows that, under some appropriate conditions, 
the eigenvalue sequence $\{\lambda_n(q)\}_{n\geqslant 1}$ 
has the desired continuity above. This result will provide a new
tool or idea for the further study of Sturm-Liouville eigenvalue problem.

The paper is structured as follows.
In Section~\ref{sec:mainthm-pre}, we present in Section~\ref{subsec:mainthm} the content of 
the main theorem, and introduce some notations in Section~\ref{subsec:ntn-pre}
as well as recalling some known facts as preliminary which are crucial for the proof of our results.
In Section~\ref{sec:prf-mainthm}, we conclude the proofs of some auxiliary lemmas, and further prove the main theorem.

\bigskip
\bigskip

\section{The main theorem and preliminary}\label{sec:mainthm-pre}
\bigskip

Throughout this paper, we  denote by $\BR$ the field of real numbers.

The symbol $L^2_{\omega}[a,b]$ denotes the weighted Hilbert space of all Lebesgue measurable, complex valued functions $f$ on $[a,b]$ satisfying $\int_a^b\omega |f|^2<\infty$ with the norm
$\norm{f}_{\omega}=(\int_a^b\omega |f|^2)^{\frac{1}{2}}$ and the inner
product $\langle f,g\rangle_\omega=\int_a^b\omega f\overline{g}$. 

We denote by $L^\infty[a,b]$ the Banach space of all essentially bounded, complex valued  functions on $[a,b]$ equipped with the canonical essential norm $\norm{\cdot}_\infty$, and by $AC[a,b]$ the space of all absolutely continuous, complex valued  functions on $[a,b]$.

\subsection{The main theorem}\label{subsec:mainthm}  ~\\

Since $p>0$ a.e. and $1/p\in L^1([a,b],\BR)$, it is easily seen that, under the following  transformation  $\CL$ of 
independent variables, called \emph{Liouville transformation} (see e.g. \cite[Page 2293]{XQ13}),
\begin{equation}\label{eq:Liouv-transf}
s=\int_a^x \frac{1}{p(t)}\, dt:= \CL(x),\quad \tilde{y}(s)=y(\CL^{-1}(s)),
\end{equation}
the problem (\ref{def:S-L}) and (\ref{def:S-L-bdry}) for $y(x)$ is rewritten as the problem for $\tilde{y}(s)$ in the form
\begin{equation}\label{def:S-L'}
-\tilde{y}''(s)+\tilde{q}(s)\tilde{y}(s)=\lambda\tilde{\omega}(s)\tilde{y}(s) \text{ on } [0,c],\quad c=\int_a^b\frac{1}{p(t)}\,dt,
\end{equation}
\begin{equation}\label{def:S-L-bdry'}
\tilde y(0)\cos\alpha+\tilde y'(0)\sin\alpha=0,\ \tilde y(c)\cos\beta+\tilde y'(c)\sin\beta=0,
\end{equation}
where $\alpha,\ \beta\in[0,\pi),$ $\tilde{q}(s)=p(\CL^{-1}(s))q(\CL^{-1}(s))$ and $\tilde{\omega}(s)=p(\CL^{-1}(s))\omega(\CL^{-1}(s))$.

It is not difficult to check that $\tilde q$, $\tilde w$ and $\tilde y$ satisfy the corresponding condition (\ref{condition}) with $[a,b]$ replaced by $[0,c]$.
More importantly, the eigenvalues of (\ref{def:S-L})-(\ref{def:S-L-bdry}) are the same as those of (\ref{def:S-L'})-(\ref{def:S-L-bdry'}).

Furthermore, for any common eigenvalue $\lambda$ of both (\ref{def:S-L})-(\ref{def:S-L-bdry}) 
and  (\ref{def:S-L'})-(\ref{def:S-L-bdry'}), denote by $E_\lambda$ and $\tilde E_\lambda$ 
the spaces of eigenfunctions associated to $\lambda$, respectively.
Then the map $y(x)\mapsto\tilde y(s)$ sets up an isometry from $E_\lambda\subset L^2_\omega [a,b]$
 onto $\tilde E_\lambda\subset L^2_{\tilde\omega}[0,c]$, and $\tilde y(s)$ on $[0,c]$ has the same range as that of $y(x)$ on $[a,b]$.

Hence, for simplicity, in the following theorem, we consider the equation (\ref{def:S-L}) for the case $p\equiv 1$ on the unit interval $[0,1]$, i.e., the eigenvalue problem
\begin{equation}\label{def:S-L''}
-y''(x)+q(x)y(x)=\lambda\omega(x)y(x) \text{ on } [0,1],
\end{equation}
\begin{equation}\label{def:S-L-bdry''}
y(0)\cos\alpha+y'(0)\sin\alpha=0,\ y(1)\cos\beta+y'(1)\sin\beta=0,
\end{equation}
where 
$$
\alpha,\ \beta\in[0,\pi),\quad   q,\ \omega\in L^1([0,1],\mathbb{R}),\quad \omega>0\ a.\ e.\text{ on }[0,1],
$$
instead of the problem (\ref{def:S-L}) and (\ref{def:S-L-bdry}). Furthermore, 
we present two hypotheses for the weight function $\omega$ of (\ref{def:S-L''})-(\ref{def:S-L-bdry''}) below:


{\bf H1:}\, $\omega(x)$ is monotonic on $[0,1]$;

{\bf H2:}\, $\inf_{x\in[0,1]} \omega(x)>0.$

\medskip
In the present paper, we mainly prove the following result.

\begin{thm}\label{thm:unif-cont-q-SL}
Suppose that the weight function $\omega$ of the eigenvalue problem (\ref{def:S-L''}) and (\ref{def:S-L-bdry''}) satisfies both of two hypothesises {\bf H1} and {\bf H2} above.
Then the eigenvalue sequence $\{\lambda_n(q)\}_{n\geqslant 1}$ of (\ref{def:S-L''})-(\ref{def:S-L-bdry''}) is uniformly locally Lipschitz continuous, in the sense of Definition~\ref{def:subunif-cont-q-SL}, with respect to the potential $q$ in $L^1([0,1],\mathbb{R})$.
\end{thm}

\bigskip

\subsection{Notations and preliminary}\label{subsec:ntn-pre}  ~\\

\medskip
For the benefit of the reader, we recall some well-known facts needed later.

\subsubsection{Differentiability of eigenvalues with respect to potential functions}\label{subsubec:deriv-q} ~\\

In this paper, by a \emph{normalized eigenfunction} of (\ref{def:S-L})-(\ref{def:S-L-bdry}) with a non-negative weight function $\omega$ we mean an eigenfunction $\varphi$ satisfying $\norm{\varphi}_\omega=1$.

\medskip
The following theorem shows the differentiability of eigenvalues of (\ref{def:S-L})-(\ref{def:S-L-bdry}) with respect to the potential functions.
\begin{thm}\label{thm:deriv} For any integer $n\geqslant 1$ and $q_0\in L^1([0,1],\mathbb{R})$, there exists a neighborhood $U(q_0)$ of $q_0$ such that, the map
$$\lambda_n: U\rightarrow\BR,\ q\in U \mapsto\lambda_n(q)\in\BR$$
is differentiable at $q_0$, and its Fr\'echet derivative is the bounded linear functional given by
\begin{equation}\label{eq:deriv}
\frac{\partial\lambda_n(q)}{\partial q}\Big\vert_{q=q_0}\cdot h=\int_0^1 \varphi_n^2(x;\lambda_n(q_0))h(x)\,dx,
\end{equation}
where $h\in L^1([0,1],\mathbb{R})$, and  $\varphi_n$ is a normalized eigenfunction associated to $\lambda_n(q_0)$ of (\ref{def:S-L})-(\ref{def:S-L-bdry}).
\end{thm}

Theorem~\ref{thm:deriv} can be viewed as a special case of a well-known theorem \cite[Theorem 4.2(6)]{KZ96} provided by Kong and Zettl.  For more details about the differentiability of eigenvalues, the reader also may refer to \cite[Theorem 3.6.1]{Zet05} and \cite{MZ96}.

\medskip
\subsubsection{Pr\"ufer transformation}\label{subsubsec:pruftransf} ~\\

Pr\"ufer transformation is an important tool in the study of Sturm-Liouville problem, 
and has several variants (see e.g. \cite{Zet05} as well as \cite{Bar57}, \cite{Reid58} and \cite{ZWMQX18}). 
In the following, we introduce the \emph{elliptic Pr\"ufer transformation}.

Consider the problem (\ref{def:S-L''}) and (\ref{def:S-L-bdry''}). Set
\begin{equation}\label{eq:pruftransf}
\rho(x;\lambda)=\sqrt{\lambda y^2(x;\lambda)+(y')^2(x;\lambda)}, \quad
\theta(x;\lambda)=\arctan\frac{\sqrt{\lambda}y(x;\lambda)}{y'(x;\lambda)},
\end{equation}
$$\lambda\geqslant 0,\quad \theta(0;\lambda)\in\BR,\quad \rho(0;\lambda)>0,$$
Then 
\begin{equation}\label{eq:pruf-sys-rad}
\theta'(x;\lambda)= \sqrt{\lambda}\left(\cos^2\theta(x;\lambda)
+\omega(x)\sin^2\theta(x;\lambda)\right)-\frac{1}{\sqrt{\lambda}}q(x)\sin^2\theta(x;\lambda)
\end{equation}
is independent of $\rho$. The equation (\ref{eq:pruf-sys-rad}) is usually called the \emph{Pr\"ufer equation}, and $\rho$ satisfies 
\begin{equation}\label{eq:pruf-sys-ang}
\rho'(x;\lambda)=\frac{\sqrt{\lambda}}{2}\left(1-\omega(x)+\frac{1}{\lambda}q(x)\right)\rho(x;\lambda)\sin 2\theta(x;\lambda).
\end{equation}

\bigskip
\bigskip

\section{The proof of Theorem~\ref{thm:unif-cont-q-SL}}\label{sec:prf-mainthm}
\bigskip

To prove our main theorem, we need to prove some lemmas and propositions as preparation.
At first, consider the initial value  problem
\begin{equation}\label{def:S-L-2}
-y''(x)=\lambda\omega(x) y(x)  \text{ on } [0,1], \ y(0)=c_1,\ y'(0)=c_2,
\end{equation}
where $c_1\geqslant 0$, $c_2\in\BR$, $\omega\in L^1([0,1],\BR)$ and $\omega>0$ a.e. on $[0,1]$.

\medskip
Applying Pr\"ufer transformation in Section~\ref{subsubsec:pruftransf} to (\ref{def:S-L-2}),
we obtain the Pr\"ufer equation for the case $q\equiv 0$ as follows:
\begin{equation}\label{eq:pruf-sys1}
\theta'(x;\lambda)=\sqrt{\lambda}(\cos^2\theta(x;\lambda)+\omega(x)\sin^2\theta(x;\lambda))
\end{equation}
with the initial condition $\theta(0;\lambda)=\arctan\frac{\sqrt{\lambda}c_1}{c_2}\in[0,\pi)$, and $\rho$ satisfies
\begin{equation}\label{eq:pruf-sys2}
\rho'(x;\lambda)=\frac{\sqrt{\lambda}}{2}\rho(x;\lambda)(1-\omega(x))\sin 2\theta(x;\lambda),\ \rho(0;\lambda)=\sqrt{\lambda c_1^2+c_2^2}\in(0,+\infty).
\end{equation}

Consequently, the solution $y(x;\lambda)$ of (\ref{def:S-L-2}) has the following expression:
\begin{equation}\label{eq:y-pruf}
y(x;\lambda)=\frac{1}{\sqrt\lambda}\rho(x;\lambda)\sin\theta(x;\lambda),
\end{equation}
where
\begin{equation}\label{eq:rho}
\rho(x;\lambda)=\rho(0;\lambda)\cdot e^{\frac{\sqrt{\lambda}}{2}\int_0^x(1-\omega(t))\sin2\theta(t;\lambda)\,dt}.
\end{equation}

Set
\begin{equation}\label{eq:H}
H(x;\lambda):=\frac{\sqrt{\lambda}}{2}\int_0^xh(t)\sin2\theta(t;\lambda)\,dt,
\end{equation}
where $h(t)=1-\omega(t)$.

\medskip

\begin{lem}\label{lem:infty-prufangl}
Let $\theta$ be defined as in (\ref{eq:pruf-sys1}) and $\omega(x)$ be in $L^1[0,1]$. If $\omega(x)\geqslant 0$ a.e. on $[0,1]$ and $\int_0^1\omega(x)\,dx>0$, then 
\begin{equation}\label{eq:infty-prufangl}
\lim_{\lambda\rightarrow +\infty}\theta(1;\lambda)=+\infty,
\end{equation}
and $\theta(x;\lambda)$ is nondecreasing  on $[0,1]$ for any fixed $\lambda\geqslant 0$.
\end{lem}

\begin{prf}
Since $\omega(x)\geqslant 0$ a.e. on $[0,1]$ and $\int_0^1\omega(x)\,dx>0$, the limit equation follows from
$$\theta(1;\lambda)-\theta(0;\lambda)=\sqrt{\lambda}\int_0^1\left(\cos^2\theta(x;\lambda)+\omega(x)\sin^2\theta(x;\lambda)\right)\,dx\geqslant\sqrt{\lambda}\int_0^1\min \{\omega(x),1\}\,dx.$$
Immediately, the remainder is proved, since the Pr\"ufer equation (\ref{eq:pruf-sys1}), together with $\omega\geqslant 0$, shows that $$\theta'(x;\lambda)\geqslant 0$$ for any $x\in [0,1]$ and $\lambda\geqslant 0$.
\end{prf}

\medskip

The following is the key lemma for the main theorem in this paper.

\begin{lem}\label{lem:bdd-RL-w}
Let $\theta$ be defined as in (\ref{eq:pruf-sys1}). Assume that both of  two hypotheses {\bf H1} and {\bf H2} hold. If $g(x): [0,1]\rightarrow\BR$ is a function whose total variation on $[0,1]$ is finite, then
\begin{equation}\label{eq:bdd-RL-w-1}
\int_0^c g(x)\sin 2\theta(x;\lambda)\, dx=O\left(\frac{1}{\sqrt{\lambda}}\right),
\end{equation}
and
\begin{equation}\label{eq:bdd-RL-w-2}
\int_0^c g(x)\cos 2\theta(x;\lambda)\, dx=O\left(\frac{1}{\sqrt{\lambda}}\right),
\end{equation}
for any $c\in [0,1]$.
\end{lem}

\begin{prf}
Here we only prove this lemma when $\omega$ is increasing. For the case that $\omega$ is decreasing, by using the transform $t=1-x$, we can keep the eigenvalues invariant, and obtain the proof in the same way.  

Since every function of bounded variation is the difference of two bounded monotonic functions,  we may further assume that $g(x)$ is monotonic.

When $c=0$, the proof is trivial. 
 
\bigskip

First, we begin to prove (\ref{eq:bdd-RL-w-1}) for $c>0$.

\smallskip

Set 
\begin{equation}\label{eq:G}
G(c;\lambda)=\frac{\sqrt\lambda}{2}\int_0^c g(x)\sin 2\theta(x;\lambda)\, dx, \quad c\in (0,1].
\end{equation}

\smallskip

{\bf Case 1}: assume that $g(x)$ is decreasing and non-negative on $[0,1]$.

\smallskip


By Lemma~\ref{lem:infty-prufangl}, for any fixed $\bar{x}\in(0,1]$ and sufficiently large $\lambda>0$ , we can find two finite sequences $$\{x_i\}_{i=0}^m,\ \{s_j\}_{j=1}^{m-1}\subseteq [0,\bar{x}],$$ such that
$$0=x_0<x_1<s_1<\cdots<x_m\leqslant\bar{x},$$ satisfying
for any $1\leqslant k\leqslant m,\ 1\leqslant j\leqslant m-1,$ $$\theta(x_k;\lambda)=k\pi,\ \theta(s_{j};\lambda)=j\pi+\frac{\pi}{2},\
m\pi\leqslant\theta(\bar x;\lambda)\leqslant (m+1)\pi,$$
and ensuring that for any $x\in(x_0,x_1),$
$$0\leqslant\theta(x;\lambda)<\pi,$$
which means that $x_1$ is the smallest one of those $x$ satisfying $\theta(x;\lambda)=\pi$.

Since $g$ is decreasing, we know that, for any integer $j\in\{1,2,...,m-1\}$,
\begin{equation}\label{ineq:mono-h-w}
g(s_j)\leqslant g(x)\leqslant g(x_j),\ \omega(x_j)\leqslant\omega(x)\leqslant\omega(s_j),\quad x\in[x_j,s_j];
\end{equation}
and
\begin{equation}\label{ineq:mono-h-w-1}
g(x_{j+1})\leqslant g(x)\leqslant g(s_{j}),\ \omega(s_{j})\leqslant\omega(x)\leqslant\omega(x_{j+1}),\quad x\in[s_{j},x_{j+1}].
\end{equation}

Moreover, by the monotonicity of $\theta$ in Lemma~\ref{lem:infty-prufangl}, we have that
$$j\pi\leqslant\theta(x;\lambda)\leqslant j\pi+\frac{\pi}{2},\quad x\in[x_{j},s_{j}];$$
and
$$j\pi+\frac{\pi}{2}\leqslant\theta(x;\lambda)\leqslant (j+1)\pi,\quad x\in[s_j,x_{j+1}].$$
Hence,
\begin{equation}\label{ineq:pos-sin2}
\sin 2\theta(x;\lambda)\geqslant 0, \quad x\in[x_j,s_j];
\end{equation}
and
\begin{equation}\label{ineq:pos-sin2-1}
\sin 2\theta(x;\lambda)\leqslant 0, \quad x\in[s_j,x_{j+1}].
\end{equation}

For simplicity, hereafter we denote $\theta(x;\lambda)$ by $\theta(x)$.

Combining the inequalities~(\ref{ineq:mono-h-w})-(\ref{ineq:pos-sin2-1}) and nonnegativity of $h$ and $w$, we obtain that
\begin{eqnarray}\label{ineq1}
\int_{x_j}^{s_j}\frac{g(s_j)\sin 2\theta(x)}{\cos^2\theta(x)+\omega(s_j)\sin^2\theta(x)}\,d\theta(x) &\leqslant& \int_{x_j}^{s_j}\frac{g(x)\sin 2\theta(x)}{\cos^2\theta(x)+\omega(x)\sin^2\theta(x)}\,d\theta(x) \nonumber \\
{} &\leqslant& \int_{x_j}^{s_j}\frac{g(x_j)\sin 2\theta(x)}{\cos^2\theta(x)+\omega(x_j)\sin^2\theta(x)}\,d\theta(x),
\end{eqnarray}
and
\begin{eqnarray}\label{ineq1-a}
\int_{s_j}^{x_{j+1}}\frac{g(s_{j})\sin 2\theta(x)}{\cos^2\theta(x)+\omega(s_{j})\sin^2\theta(x)}\,d\theta(x) &\leqslant& \int_{s_{j}}^{x_{j+1}}\frac{g(x)\sin 2\theta(x)}{\cos^2\theta(x)+\omega(x)\sin^2\theta(x)}\,d\theta(x) \nonumber \\
{} &\leqslant& \int_{s_{j}}^{x_{j+1}}\frac{g(x_{j+1})\sin 2\theta(x)}{\cos^2\theta(x)+\omega(x_{j+1})\sin^2\theta(x)}\,d\theta(x).
\end{eqnarray}

Define an auxiliary function as follows:
\begin{equation}\label{def:funct}
f(t)=\int_0^{\frac{\pi}{2}}\frac{g(t)\sin 2u}{\cos^2u+\omega(t)\sin^2u}\,du,\quad t\in[0,1].
\end{equation}

Then, substituting $\theta(x)$ for $u$, by the periodicity of $\sin 2u$, we have that
\begin{equation}\label{def:funct}
f(t)=\int_{x_j}^{s_j}\frac{g(t)\sin 2\theta(x)}{\cos^2\theta(x)+\omega(t)\sin^2\theta(x)}\,d\theta(x)=-\int_{s_j}^{x_{j+1}}\frac{g(t)\sin 2\theta(x)}{\cos^2\theta(x)+\omega(t)\sin^2\theta(x)}\,d\theta(x).
\end{equation}

So, it follows from (\ref{ineq1})-(\ref{ineq1-a}) that
\begin{equation}\label{ineq1-1}
f(s_j)\leqslant
\int_{x_j}^{s_j}\frac{g(x)\sin 2\theta(x)}{\cos^2\theta(x)+\omega(x)\sin^2\theta(x)}\,d\theta(x)\leqslant
f(x_j)
\end{equation}
and
\begin{equation}\label{ineq2}
-f(s_j)\leqslant
\int_{s_j}^{x_{j+1}}\frac{g(x)\sin 2\theta(x)}{\cos^2\theta(x)+\omega(x)\sin^2\theta(x)}\,d\theta(x)\leqslant
-f(x_{j+1}),
\end{equation}
where $1\leqslant j\leqslant m-1.$

Adding together the two inequalities above, we have that
\begin{equation}\label{ineq12}
0\leqslant
\int_{x_{j}}^{x_{j+1}}\frac{g(x)\sin 2\theta(x)}{\cos^2\theta(x)+\omega(x)\sin^2\theta(x)}\,d\theta(x)\leqslant
f(x_{j})-f(x_{j+1}), 
\end{equation}
where $1\leqslant j\leqslant m-1$.

For the last interval $[x_m,\bar{x}]$, it can be known, from the similar argument as above, that
\begin{equation}\label{ineq3}
\left\{ \begin{array}{lll}
f(\bar{x})\leqslant 
\int_{x_m}^{\bar{x}}\frac{g(x)\sin 2\theta(x)}{\cos^2\theta(x)+\omega(x)\sin^2\theta(x)}\,d\theta(x)\leqslant f(x_m),   &    \textrm{if  $m\pi\leqslant\theta(\bar{x})\leqslant m\pi+\frac{\pi}{2}$};\\

{} & {}\\

0\leqslant 
\int_{x_m}^{\bar{x}}\frac{g(x)\sin 2\theta(x)}{\cos^2\theta(x)+\omega(x)\sin^2\theta(x)}\,d\theta(x)\leqslant f(x_m)-f(\bar{x}),       &   \textrm{if $m\pi+\frac{\pi}{2}<\theta(\bar{x})\leqslant(m+1)\pi$}.
\end{array}  \right.
\end{equation}

From monotonicity and non-negativity of $h$ and $\omega$, it is apparent that $f(t)$ is non-negative and decreasing on $[0,1]$, and so
\begin{equation}\label{ineq:f}
0\leqslant f(t)\leqslant f(0)=\int_0^{\frac{\pi}{2}}\frac{g(0)\sin 2u}{\cos^2u+\omega(0)\sin^2u}\,du<+\infty,
\end{equation}
where the finiteness of the integral in (\ref{ineq:f}) owes to $\omega(0)>0$.

Then, it follows from (\ref{ineq12})--(\ref{ineq3}) that, for the $\bar{x}$ arbitrarily given above,
\begin{equation}\label{ineq4-1}
0 \leqslant  \left(\sum_{j=1}^{m-1}\int_{x_j}^{x_{j+1}}+\int_{x_m}^{\bar{x}}\right)\frac{g(x)\sin 2\theta(x)}{\cos^2\theta(x)+\omega(x)\sin^2\theta(x)}\,d\theta(x) \leqslant  f(x_1)\leqslant  f(0).
\end{equation}

Moreover, since $\omega(0)>0$,  we also have that

\begin{equation}\label{ineq4-2}
\abs{\int_{x_0}^{x_1}\frac{g(x)\sin 2\theta(x)}{\cos^2\theta(x)+\omega(x)\sin^2\theta(x)}\,d\theta(x)}
\leqslant\frac{\pi g(0)}{\min\{\omega(0),1\}}.
\end{equation}

Notice that
\begin{eqnarray}\label{eq:G-decomp}
G(\bar x;\lambda) &=& \frac{\sqrt{\lambda}}{2}\int_0^{\bar{x}}g(x)\sin2\theta(x;\lambda)\,dt \nonumber \\
{} &=& \frac{\sqrt{\lambda}}{2}\left(\int_{x_0}^{x_1}+\sum_{j=1}^{m-1}\int_{x_{j}}^{x_{j+1}}+\int_{x_m}^{\bar{x}}\right)g(x)\sin 2\theta(x)\,dx  \nonumber \\
{} &=& \frac{1}{2}\left(\int_{x_0}^{x_1}+\sum_{j=1}^{m-1}\int_{x_{j}}^{x_{j+1}}+\int_{x_m}^{\bar{x}}\right)\frac{g(x)\sin 2\theta(x)}{\cos^2\theta(x)+\omega(x)\sin^2\theta(x)}\,d\theta(x)
\end{eqnarray}

Set $$G_0:=f(0)+\frac{\pi g(0)}{\min\{\omega(0),1\}}.$$

Therefore, from (\ref{ineq4-1})--(\ref{eq:G-decomp}) and the arbitrariness of $\bar x$, we can derive that, for any $c\in[0,1]$ and sufficiently large $\lambda>0$, one has
\begin{equation}\label{ineq:G_0}
\abs{G(c;\lambda)}\leqslant \frac{G_0}{2}<+\infty,
\end{equation}
which implies (\ref{eq:bdd-RL-w-1}) in  {\bf Case 1}.

\smallskip

{\bf Case 2}: assume that $g(x)$ is decreasing  on $[0,1]$, but is not needed to be non-negative.

\smallskip

Let $u(x)=g(x)-g(1)$ and $v(x)\equiv 1$. Then $u(x)$ is non-negative and also decreasing on $[0,1]$, and $g(x)=u(x)+g(1)v(x)$. So, for any $c\in [0,1]$, we have
$$\int_0^c g(x)\sin 2\theta(x;\lambda)\, dx=\int_0^c u(x)\sin 2\theta(x;\lambda)\, dx+g(1)\int_0^c  v(x)\sin 2\theta(x;\lambda)\, dx.$$

Applying the result in {\bf Case 1} to  the functions $u$ and $v$, we obtain (\ref{eq:bdd-RL-w-1}) in {\bf Case 2}.

\smallskip

{\bf Case 3}: assume that $g(x)$ is increasing  on $[0,1]$. 

\smallskip
Set $u(x)=g(1)-g(x)$ and $v(x)\equiv 1$.  
So $g(x)=g(1) v(x)-u(x)$, and $u(x)$ is decreasing on $[0,1]$. Then (\ref{eq:bdd-RL-w-1}) follows from the trick similar to that in {\bf Case 2}.

\smallskip
From the argument above, the proof of (\ref{eq:bdd-RL-w-1}) is done.

\bigskip
For (\ref{eq:bdd-RL-w-2}), the proof is similar to that of (\ref{eq:bdd-RL-w-1}).

\smallskip
Set 
\begin{equation}\label{eq:G1}
\tilde G(c;\lambda)=\frac{\sqrt\lambda}{2}\int_0^c g(x)\cos 2\theta(x;\lambda)\, dx, \quad c\in (0,1].
\end{equation}

By Lemma~\ref{lem:infty-prufangl}, for any fixed $\bar{x}\in(0,1]$ and sufficiently large $\lambda>0$ , we also can find two finite sequences $$\{x_i\}_{i=0}^m,\ \{s_j\}_{j=1}^{m-1}\subseteq [0,\bar{x}],$$ such that
$$0=x_0<x_1<s_1<\cdots<x_m\leqslant\bar{x},$$ satisfying
for any $1\leqslant k\leqslant m,\ 1\leqslant j\leqslant m-1,$ $$\theta(x_k;\lambda)=k\pi+\frac{\pi}{4},\ \theta(s_{j};\lambda)=j\pi+\frac{3\pi}{4},\
m\pi+\frac{\pi}{4}\leqslant\theta(\bar x;\lambda)\leqslant (m+1)\pi+\frac{\pi}{4},$$
and ensuring that for any $x\in(x_0,x_1),$
$$0\leqslant\theta(x;\lambda)<\frac{5\pi}{4},$$
which means that $x_1$ is the smallest one of those $x$ satisfying $\theta(x;\lambda)=\frac{5\pi}{4}$.

Then
\begin{equation}\label{eq:G1-decomp}
\tilde G(\bar x;\lambda) = \frac{1}{2}\left(\int_{x_0}^{x_1}+\sum_{j=1}^{m-1}\int_{x_{j}}^{x_{j+1}}+\int_{x_m}^{\bar{x}}\right)\frac{g(x)\cos 2\theta(x)}{\cos^2\theta(x)+\omega(x)\sin^2\theta(x)}\,d\theta(x).
\end{equation}

Moreover, it can be clearly seen that
$$j\pi+\frac{\pi}{4}\leqslant\theta(x;\lambda)\leqslant j\pi+\frac{3\pi}{4},\quad x\in[x_{j},s_{j}];$$
and
$$j\pi+\frac{3\pi}{4}\leqslant\theta(x;\lambda)\leqslant (j+1)\pi+\frac{\pi}{4},\quad x\in[s_j,x_{j+1}].$$
Hence,
\begin{equation}\label{ineq:pos-cos2}
\cos 2\theta(x;\lambda)\leqslant 0, \quad x\in[x_j,s_j];
\end{equation}
and
\begin{equation}\label{ineq:pos-cos2-1}
\cos 2\theta(x;\lambda)\geqslant 0, \quad x\in[s_j,x_{j+1}].
\end{equation}

Similar to (\ref{def:funct}), we define the corresponding auxiliary function as follows:
\begin{equation}\label{def:funct1}
\tilde f(t)=\int_{\frac{\pi}{4}}^{\frac{3\pi}{4}}\frac{g(t)\cos 2u}{\cos^2u+\omega(t)\sin^2u}\,du,\ \quad t\in[0,1].
\end{equation}
Note that $f(t)\in L^\infty[0,1]$.

\smallskip

Under the above setting (\ref{eq:G1})-(\ref{def:funct1}), from lines of argument similar to those of (\ref{eq:bdd-RL-w-1}), it can be shown that, when $g(x)$ is decreasing and non-negative on $[0,1]$, for any $c\in [0,1]$ and sufficiently large $\lambda>0$, one has

\begin{equation}\label{eq:G1-0}
\abs{\tilde G(c;\lambda)}\leqslant\frac{1}{2}\left (\abs{f(0)}+\frac{5\pi}{4}\cdot\frac{g(0)}{\min\{\omega(0),1\}}\right),
\end{equation}
which implies (\ref{eq:bdd-RL-w-2}) in {\bf Case 1}.

Similarly, (\ref{eq:bdd-RL-w-2}) in both of {\bf Case 2} and {\bf Case 3} also can be obtained.

This lemma is proved.

\end{prf}

\medskip

The following result is a direct consequence of Lemma~\ref{lem:bdd-RL-w}.

\begin{lem}\label{lem:bdd-H}
If both of {\bf H1} and {\bf H2} hold, then  $H(x;\lambda)$ in (\ref{eq:H}) is uniformly bounded for all sufficiently large $\lambda>0$.
\end{lem}

\medskip

The next lemma can be considered as an analogue of Riemann-Lesbegue lemma.

\begin{lem}\label{lem:RL-w}
Let $\theta$ be defined as in (\ref{eq:pruf-sys1}). Assume that {\bf H1} holds. If $g(x)$ is an arbitrary element in $L^1([0,1],\BR)$, then
\begin{equation}\label{eq:RL-w-1}
\lim_{\lambda\rightarrow +\infty}\int_0^c g(x)\sin 2\theta(x;\lambda)\, dx=0,
\end{equation}
and
\begin{equation}\label{eq:RL-w-2}
\lim_{\lambda\rightarrow +\infty}\int_0^c g(x)\cos 2\theta(x;\lambda)\, dx=0,
\end{equation}
for any $c\in [0,1]$.
\end{lem}

\begin{prf} 
First,  we claim that both of (\ref{eq:RL-w-1}) and (\ref{eq:RL-w-2}) hold when $g(x)\in AC[0,1]$ and $\omega$ satisfies both of two hypotheses {\bf H1} and {\bf H2}.  Since every absolutely continuous function has bounded variation, this claim is obviously true because of  Lemma~\ref{lem:bdd-RL-w}.

Next, we retain {\bf H1} but  remove  {\bf H2}.  Set $\omega_n=\omega+\frac{1}{n}$.   
Since all of absolutely continuous functions are dense in $L^1[0,1]$,  we can deduce by the above claim that (\ref{eq:RL-w-1})-(\ref{eq:RL-w-2}) hold for every $\omega_n$ and any $g(x)\in L^1[0,1]$. 
Consequently, (\ref{eq:RL-w-1})-(\ref{eq:RL-w-2}) are true for any non-negative monotonic weight $\omega$ and integrable function $g(x)$, since $\omega_n$ is uniformly convergent to $\omega$ on $[0,1]$ and the Pr\"ufer argument $\theta$ depends continuously on its weight function (see \cite[Theorem~4.5.1]{Zet05}).
\end{prf}


\medskip

By the above lemmas, we can establish, on any bounded subset of $L^1([0,1],\BR)$, the uniform boundedness of the normalized eigenfunctions of the eigenvalue problem (\ref{def:S-L''})-(\ref{def:S-L-bdry''}).

\begin{prop}\label{prop:SL-solution-bdd} Consider the eigenvalue problem (\ref{def:S-L''})-(\ref{def:S-L-bdry''}), and suppose that the weight function $\omega$ satisfies {\bf H1}--{\bf H2}.
Then, for any $L^1$-norm bounded subset $\Omega$ of $L^1([0,1],\BR)$, there exists a positive number $M(\Omega)$ such that, for any normalized eigenfunction $\varphi_n(x;\lambda_n(q))$ of  (\ref{def:S-L''})-(\ref{def:S-L-bdry''}), one has $$\abs{\varphi_n(x;\lambda_n(q))}\leqslant M(\Omega)$$ for all $n\geqslant 1$, $q\in\Omega$ and $x\in[0,1]$.
\end{prop}

\begin{prf}
Consider the initial value problem as follows:
\begin{equation}\label{def:S-L-1}
-y''(x)+q(x)y(x)=\lambda\omega(x) y(x) \text{ on } [0,1],\ 
y(0)=C_1,\ y'(0)=C_2,
\end{equation}
where $q,\ \omega \in L^1([0,1],\BR)$, $\omega>0$ a.e. on $[0,1]$ and $C_1,\ C_2$ are two arbitrary fixed real numbers satisfying 
\begin{equation}\label{def:S-L-bdry-1'}
C_1\cos\alpha+C_2\sin\alpha=0,\quad y(1)\cos\beta+y'(1)\sin\beta=0
\end{equation}
where $\alpha$ and $\beta$ are given in the boundary condition (\ref{def:S-L-bdry''}).

We may as well assume that $C_1\neq 0$.

Choose two linearly independent solutions $\phi$ and $\psi$ of (\ref{def:S-L-2}), such that
$\phi(0)=0,\ \phi'(0)=1$ and $\psi(0)=1,\ \psi'(0)=0$. Clearly, Wronskian determinant $$W[\psi,\phi]=\det\begin{pmatrix}
   \psi(x)   &  \phi(x) \\
   \psi'(x)   &  \phi'(x)
\end{pmatrix}$$ of $\psi$ and $\phi$ equals to 1.

We may choose $\psi$ and $\phi$ as follows:
\begin{equation}\label{def:solution-phi-psi}
\left\{ \begin{array}{ll}
\phi(x;\lambda)=\frac{1}{\sqrt{\lambda}}r(x;\lambda)\sin\nu(x;\lambda),   &    r(0;\lambda)=1,\ \nu(0;\lambda)=0;\\
\psi(x;\lambda)=\mu(x;\lambda)\sin\sigma(x;\lambda),    &    \mu(0;\lambda)=1,\ \sigma(0;\lambda)=\frac{\pi}{2}.
\end{array}  \right.
\end{equation}
where $(r,\nu)$ and $(\mu,\sigma)$ satisfies the corresponding equation (\ref{eq:rho}).

So, by Pr\"ufer transformation, we obtain that
\begin{equation}\label{eq:solution-phi-psi}
\left\{ \begin{array}{ll}
\phi'(x;\lambda)=r(x;\lambda)\cos\nu(x;\lambda),   &   \phi(0;\lambda)=0,\ \phi'(0;\lambda)=1;\\
\psi'(x;\lambda)=\sqrt{\lambda}\mu(x;\lambda)\cos\sigma(x;\lambda),    &  \psi(0;\lambda)=1,\ \psi'(0;\lambda)=0.
\end{array}  \right.
\end{equation}

For the initial condition in (\ref{def:S-L-1}), using the formula of variation of constant, we can derive that the unique solution $y(x;\lambda)$ of (\ref{def:S-L-1}) satisfies the integral equation 
\begin{equation}\label{eq:var-const}
y(x;\lambda)=C_1\psi(x;\lambda)+C_2\phi(x;\lambda)+\int_0^x[\phi(x;\lambda)\psi(t;\lambda)-\phi(t;\lambda)\psi(x;\lambda)]q(t)y(t;\lambda)\,dt.
\end{equation}

Putting  (\ref{def:solution-phi-psi})  into (\ref{eq:var-const}), we have
\begin{equation}\label{eq:var-const-polar}
y(x;\lambda)=C_1\mu(x;\lambda)\sin\sigma(x;\lambda)+C_2\frac{1}{\sqrt{\lambda}}r(x;\lambda)\sin\nu(x;\lambda)+\frac{1}{\sqrt{\lambda}}\int_0^xR(x;t;\lambda)q(t)y(t;\lambda)\,dt,
\end{equation}
where $R(x;t;\lambda)=r(x;\lambda)\sin\nu(x;\lambda)\mu(t;\lambda)\sin\sigma(t;\lambda)-r(t;\lambda)\sin\nu(t;\lambda)\mu(x;\lambda)\sin\sigma(x;\lambda)$.

Because $(r,\nu)$ and $(\mu,\sigma)$ satisfy the corresponding equation (\ref{eq:rho}), it is easily known from Lemma~\ref{lem:bdd-H} that, there exists positive numbers $M_0$ and $K$ such that, for any $\lambda\geqslant K$, 
\begin{equation}\label{ineq:rhotheta}
\abs{r(x;\lambda)\sin\nu(x;\lambda)}\leqslant M_0,
\end{equation}
and
\begin{equation}\label{ineq:musigma}
\abs{\mu(x;\lambda)\sin\sigma(x;\lambda)}\leqslant M_0,
\end{equation}
and then, 
\begin{equation}\label{ineq:R}
\abs{R(x;t;\lambda)}\leqslant 2M_0^2.
\end{equation}

Set $$B(\Omega):=\sup\{\norm{q}_{L^1} |\ q\in\Omega\}.$$ 
By the inequalities (\ref{eq:var-const-polar})-(\ref{ineq:R})and Gronwall inequality (see e.g. \cite[Theorem 1.4.1(i)]{Zet05} and \cite[Theorem 1.3.2]{Pach98}), it is apparent that, for any $\lambda\geqslant \max\{1, K\}$,
\begin{eqnarray}\label{ineq:y-bdd}
\abs{y(x;\lambda)} &\leqslant& \abs{C_1\mu\sin\sigma}+\frac{\abs{C_2r\sin\nu}}{\sqrt{\lambda}}+\int_0^x\frac{\abs{R(x;t;\lambda)q(t)}}{\sqrt{\lambda}}y(t;\lambda)\,dt,  \nonumber \\
{} &\leqslant& (\abs{C_1}+\abs{C_2})M_0+\int_0^x \abs{R(x;t;\lambda)q(t)}y(t;\lambda)\,dt  \nonumber \\
{} &\leqslant& (\abs{C_1}+\abs{C_2})M_0\cdot e^{\int_0^x \abs{R(x;s;\lambda)q(s)}\,ds} \nonumber \\
{} &\leqslant& (\abs{C_1}+\abs{C_2})M_0\cdot e^{2M_0^2B(\Omega)},
\end{eqnarray}
which implies that
\begin{equation}\label{eq:Rqy-bdd}
\abs{\int_0^xR(x;t;\lambda)q(t)y(t;\lambda)\,dt}\leqslant (\abs{C_1}+\abs{C_2})M_0(e^{2M_0^2B(\Omega)}-1).
\end{equation}

By (\ref{eq:var-const-polar}), (\ref{ineq:rhotheta}) and (\ref{eq:Rqy-bdd}),  it can be seen that, for any $\lambda\geqslant \max\{1, K\}$,
\begin{equation}\label{eq:y-asym}
y(x;\lambda)=C_1\mu(x;\lambda)\sin\sigma(x;\lambda)+O(\frac{1}{\sqrt{\lambda}}).
\end{equation}

And then,  there exists a positive number $M_1$ such that, for any $\lambda\geqslant \max\{1, K\}$,
\begin{equation}\label{ineq:y-sq-asym}
C_1^2\mu^2(x;\lambda)\sin^2\sigma(x;\lambda)-\abs{y(x;\lambda)}^2\leqslant  \frac{M_1}{\sqrt{\lambda}}.
\end{equation}

Let $\{\lambda_n\}_{n\geqslant 1}$ be the eigenvalue sequence of the eigenvalue problem (\ref{def:S-L''})-(\ref{def:S-L-bdry''}). Then the unique solution $y(x;\lambda_n)$ of the initial value problem (\ref{def:S-L-1}) is also a eigenfunction of (\ref{def:S-L''})-(\ref{def:S-L-bdry''}) corresponding to $\lambda_n$. So we can find a number $\beta(\lambda_n)$ such that $$\varphi_n(x;\lambda_n)=\beta(\lambda_n)y(x;\lambda_n)$$
satisfying $$\int_0^1\omega(x)\abs{\varphi_n(x;\lambda_n)}^2\,dx=1. $$

Thereupon, we have 
\begin{equation}\label{eq:beta-y}
\beta^2(\lambda_n)\int_0^1\omega(x)\abs{y(x;\lambda_n)}^2\,dx=1.
\end{equation}

Since $\lambda_n\rightarrow +\infty$ as $n\rightarrow +\infty$, there exists a sufficiently large positive integer $N_0$ such that $\lambda_n\geqslant\max\{1, K\}$ for any $n\geqslant N_0$. 

Hence, by (\ref{ineq:y-sq-asym}) and (\ref{eq:beta-y}), we have
\begin{equation}\label{eq:y-beta-bdd}
\beta^2(\lambda_n)\int_0^1\omega(x)\mu^2(x;\lambda_n)\sin^2\sigma(x;\lambda_n)\,dx\leqslant \frac{1}{C_1^2}+\frac{M_1\norm{\omega}_{L^1}}{C_1^2\sqrt{\lambda_n}}\beta^2(\lambda_n),
\end{equation}
for any $n\geqslant N_0$.

Since $(\mu,\sigma)$ satisfies the corresponding equation (\ref{eq:rho}), the equation (\ref{eq:y-beta-bdd}), together with $\mu(0;\lambda)=1$, yields that,
\begin{equation}\label{eq:alpha}
\beta^2(\lambda_n)\int_0^1\omega(x)e^{2H(x;\lambda_n)}\sin^2\sigma(x;\lambda_n)\,dx\leqslant \frac{1}{C_1^2}+\frac{M_1\norm{\omega}_{L^1}}{C_1^2\sqrt{\lambda_n}}\beta^2(\lambda_n),
\end{equation}
for any $n\geqslant N_0$.

Lemma~\ref{lem:bdd-H} tells us that, there exists a positive number $H_0$ and a sufficiently large $N_1\geqslant N_0$ such that
\begin{equation}\label{eq:est_hsin}
e^{-\frac{H_0}{2}}\leqslant e^{H(x;\lambda_n)}\leqslant e^{\frac{H_0}{2}},\text{ for any $x\in[0,1]$ and $n\geqslant N_1$}.
\end{equation}

Then, for any $n\geqslant N_1$, we have
\begin{equation}\label{ineq:alpha}
\beta^2(\lambda_n)\int_0^1\omega(x)\sin^2\sigma(x;\lambda_n)\,dx\leqslant \frac{e^{H_0}}{C_1^2}+\frac{e^{H_0}M_1\norm{\omega}_{L^1}}{C_1^2\sqrt{\lambda_n}}\beta^2(\lambda_n),
\end{equation}
that is,
\begin{equation}\label{ineq:alpha1}
\beta^2(\lambda_n)\int_0^1\omega(x)\frac{1-\cos 2\sigma(x;\lambda_n)}{2}\,dx\leqslant  \frac{e^{H_0}}{C_1^2}+\frac{e^{H_0}M_1\norm{\omega}_{L^1}}{C_1^2\sqrt{\lambda_n}}\beta^2(\lambda_n).
\end{equation}

So, by Lemma~\ref{lem:RL-w} and (\ref{eigen-infty}), for any fixed $\gamma\in(0,1)$, we can choose a sufficiently large integer $N\geqslant N_1$, such that, as long as $n\geqslant N$,  one has
\begin{equation}\label{ineq:RL}
\int_0^1\omega(x)\cos 2\sigma(x;\lambda_n)\, dx\leqslant\gamma\int_0^1\omega
\end{equation}
and
\begin{equation}\label{ineq:lambda_n}
\frac{e^{H_0}M_1}{C_1^2\sqrt{\lambda_n}}< \frac{(1-\gamma)}{4}.
\end{equation}

Consequently, by the inequalities~(\ref{ineq:alpha1})-(\ref{ineq:lambda_n}), we have
\begin{equation}\label{ineq:beta-bdd}
\abs{\beta(\lambda_n)}<\frac{1}{\abs{C_1}}\frac{2e^{\frac{H_0}{2}}}{\sqrt{(1-\gamma)\norm{\omega}_{L^1}}}.
\end{equation}

Set $$M(\Omega):=\max\left\{M_0e^{2M_0^2B(\Omega)}\left(1+\abs{\frac{C_2}{C_1}}\right)\frac{2e^{\frac{H_0}{2}}}{\sqrt{(1-\gamma)\norm{\omega}_{L^1}}},\ \norm{\varphi_1}_\infty, ... ,\norm{\varphi_{N-1}}_\infty\right\},$$  where $\varphi_i$ is the unique normalized eigenfunction corresponding to the $i$th eigenvalue $\lambda_i$.

Hence, it follows from (\ref{ineq:y-bdd}) and (\ref{ineq:beta-bdd}) that, for any $n\geqslant 1$, 
$$\abs{\varphi_n(x;\lambda_n)}=\abs{\beta(\lambda_n)y(x;\lambda_n)}\leqslant M(\Omega).$$ 

The proof is finished.
\end{prf}

\medskip
Now, it's time to give the proof of Theorem~\ref{thm:unif-cont-q-SL}.
\bigskip

\noindent \textbf{Proof of Theorem~\ref{thm:unif-cont-q-SL}:} 
Set $\bar B_M:=\{q\in L^1[0,1]\ |\ \norm{q}_{L^1}\leqslant M\}$, which is convex.
For any $L^1$-norm bounded subset $\Omega$ of $L^1([0,1],\BR)$, set $B(\Omega):=\sup\{\norm{q}_{L^1} |\ q\in\Omega\}$. 
It is easily seen that $\Omega\subset \bar B_{B(\Omega)}$. Hence we only need to prove our result holds for convex sets.

Let $\Omega$ be an arbitrary convex $L^1$-norm bounded subset of $L^1([0,1],\BR)$. For any two $q_1,\ q_2\in\Omega$ and $\Delta q=q_2-q_1$, set
$$q_t(x)=q_1(x)+t\cdot\Delta q(x)$$ and $$\tilde\lambda_n(t)=\lambda_n(q_t),\ t\in[0,1].$$

Let $\varphi_n(x;t)$ be the unique normalized eigenfunction of $\tilde\lambda_n(t)$.
By Theorem~\ref{thm:deriv}, it is apparent that
\begin{equation}\label{equ:deriv1}
\frac{\partial\lambda_n(q_t)}{\partial q_t}=\varphi_n^2(x;t)
\end{equation}
as a bounded linear functional on $L^1([0,1],\BR)$.

Then, by (\ref{equ:deriv1}),  we obtain that
\begin{eqnarray}\label{ineq:main}
\abs{\lambda_n(q_2)-\lambda_n(q_1)} & = & \abs{\tilde\lambda_n(1)-\tilde\lambda_n(0)}=\abs{\int_0^1\frac{d\lambda_n(q_t)}{dt}\,dt}  \nonumber \\
{} & = & \abs{\int_0^1\frac{\partial\lambda_n(q)}{\partial q}\Big\vert_{q=q_t}\cdot\frac{d(q_t)}{dt}\,dt} \nonumber \\
{} &=& \abs{\int_0^1\frac{\partial\lambda_n(q)}{\partial q}\Big\vert_{q=q_t}\cdot\Delta q(x)\,dt}  \nonumber \\
{} &=& \abs{\int_0^1\int_0^1 \varphi_n^2(x;t)\Delta q(x)\,dxdt}  \nonumber \\
{} &\leqslant& \int_0^1\int_0^1\varphi_n^2(x;t)\abs{\Delta q(x)}\,dxdt.
\end{eqnarray}

Finally, due to Proposition~\ref{prop:SL-solution-bdd} and (\ref{ineq:main}), the proof is done.
\hfill $\Box$\\ \smallskip

\bigskip

\section*{Acknowledgement}

\medskip

The authors gratefully acknowledge the anonymous referees for their valuable comments which substantially improved the quality of this paper. It is especially helpful that a vital error about the proof of Lemma~\ref{lem:RL-w} in the original manuscript was point out by the referees.
The authors also would like to thank Professor Bing Xie (Shandong University, Weihai) and Dr. Qianhong Huang (University of Alberta, Canada) for some helpful discussions and  suggestions.

\end{document}